\documentclass[11pt,oneside]{article}%
\usepackage{amsmath}
\usepackage{amsfonts}
\usepackage{amssymb}
\usepackage{graphicx}
\usepackage[dvipsnames]{xcolor}
\usepackage[numbers]{natbib}
\usepackage[left=1in,right=1in,top=1.1in,bottom=1.1in]{geometry}
\usepackage[colorlinks=true,linkcolor=RawSienna,citecolor=RawSienna,urlcolor=RawSienna,bookmarksopen=true,pdfstartview=FitB]{hyperref}
\usepackage[onehalfspacing]{setspace}

\newtheorem{theorem}{Theorem}

\newtheorem{remark}{Remark}

\allowdisplaybreaks
\begin{document}

\title{\vspace{-1.2cm}A strong law of large numbers for simultaneously testing parameters of
Lancaster bivariate distributions}
\author{Xiongzhi Chen\thanks{Department of Mathematics and Statistics, Washington
State University, Pullman, WA 99164, USA. Email:
\texttt{xiongzhi.chen@wsu.edu}.} }
\date{}
\maketitle

\begin{abstract}
We prove a strong law of large numbers for simultaneously testing parameters
of a large number of dependent, Lancaster bivariate random variables with
infinite supports, and discuss its implications.

\end{abstract}

\section{Introduction}

Multiple hypothesis testing with \textit{false discovery rate} (FDR,
\cite{Benjamini:1995}) control has been widely applied to various scientific
endeavors, and it can often be stated as follows. There are $m\in\mathbb{N}$
test statistics $\left\{  \zeta_{i}\right\}  _{i=1}^{m}$ such that $\zeta_{i}$
has parameter $\mu_{i}$, and the $i$th null hypothesis is $H_{i0}:\mu_{i}%
=\mu_{0}$ (versus its alternative hypothesis $H_{i1}:\mu_{i}\neq\mu_{0}$) for
a fixed, known $\mu_{0}\in\Theta\subseteq\mathbb{R}$, where $\Theta$ is the
parameter space for the $m$ $\zeta_{i}$'s. Define $p_{i}=1-F_{i}\left(
\zeta_{i}\right)  $ as the one-sided p-value for $\zeta_{i}$, where $F_{i}$ is
the \textit{cumulative distribution function} (CDF) of $\zeta_{i}$ when
$\mu_{i}=\mu_{0}$. Let $I_{0,m}$ be the set of indices of the true null
hypotheses, and denote its cardinality (often being positive) by $m_{0}$.
Consider the multiple testing procedure (MTP) with a fixed rejection threshold
$t\in\left[  0,1\right]  $ that rejects $H_{i0}$ \textit{if and only if} (iff)
$p_{i}\leq t$. Then the MTP induces $R_{m}\left(  t\right)  =\sum_{i=1}%
^{m}1\left\{  p_{i}\leq t\right\}  $, the number of rejections, and
$V_{m}\left(  t\right)  =\sum_{i\in I_{0,m}}1\left\{  p_{i}\leq t\right\}  $,
the number of false discoveries, where $1A$ is the indicator function of a set
$A$. Further, the \textit{false discovery proportion} (FDP) and FDR of the MTP
are
\begin{equation}
\mathrm{FDP}_{m}\left(  t\right)  =\frac{V_{m}\left(  t\right)  }{R_{m}\left(
t\right)  \vee1}\text{ \ \ \ and \ \ \ \ }\mathrm{FDR}_{m}\left(  t\right)
=\mathbb{E}\left[  \mathrm{FDP}_{m}\left(  t\right)  \right]  \label{eqFDR}%
\end{equation}
respectively, where the operator $\vee$ returns the maximum of its two
arguments. When $m$, the number of tests to conduct, is large, we aim to
control the FDR of the MTP at a given level $\theta\in\left(  0,1\right)  $ by
choosing an appropriate $t$ or to estimate the FDP or FDR of the MTP at a
given threshold $t$.

However, the test statistics $\left\{  \zeta_{i}\right\}  _{i=1}^{m}$ are
often dependent on each other, and under dependence the behavior of the FDP is
usually unstable and can sometimes be unpredictable; see, e.g.,
\cite{Finner:2007}, \cite{Owen:2005} and \cite{Schwartzman2011}. This can make
irreproducible and untrustable the inferential results from the MTP. The very
few works of \cite{Schwartzman:2015}, \cite{Chen:2014SLLN},
\cite{Delattre:2016} and \cite{Fan:2012} studied the asymptotic behavior of
$R_{m}\left(  t\right)  $ or $m^{-1}R_{m}\left(  t\right)  $ under dependence
by utilizing conditions on the correlation matrix $\mathbf{R}=\left(
\rho_{ij}\right)  $ of $\boldsymbol{\zeta}=\left(  \zeta_{1},\ldots,\zeta
_{m}\right)  $. However, they all considered the setting where each $\zeta
_{i}$ is a Gaussian random variable. Specifically, when each pair $\left(
\zeta_{i},\zeta_{j}\right)  ,i\neq j$ is bivariate Gaussian, the authors of
\cite{Chen:2014SLLN} proved \textquotedblleft a SLLN for $R_{m}\left(
t\right)  $ and $V_{m}\left(  t\right)  $\textquotedblright, i.e.,

\begin{description}
\item[C1)] If%
\begin{equation}
m^{-2}\left\Vert \mathbf{R}\right\Vert _{1}=O\big(m^{-\delta}\big)\text{ \ for
some \ }\delta>0,\label{eq1}%
\end{equation}
then%
\begin{equation}
\left\{
\begin{array}
[c]{c}%
\lim_{m\rightarrow\infty}\left\vert m^{-1}R_{m}\left(  t\right)
-\mathbb{E}\left[  m^{-1}R_{m}\left(  t\right)  \right]  \right\vert =0\text{
\ almost surely,}\\
\lim_{m\rightarrow\infty}\left\vert m^{-1}V_{m}\left(  t\right)
-\mathbb{E}\left[  m^{-1}V_{m}\left(  t\right)  \right]  \right\vert =0\text{
\ almost surely.}%
\end{array}
\right.  \label{eqSLLNR}%
\end{equation}

\item[C2)] If $\liminf_{m\rightarrow\infty}m_{0}m^{-1}>0$ and (\ref{eq1})
hold, then%
\begin{equation}
\lim_{m\rightarrow\infty}\left\vert m_{0}^{-1}V_{m}\left(  t\right)
-\mathbb{E}\left[  m_{0}^{-1}V_{m}\left(  t\right)  \right]  \right\vert
=0\text{ almost surely.}\label{eqSLLNV}%
\end{equation}

\end{description}

\noindent Here \textquotedblleft the $l_{1}$-norm $\left\Vert \mathbf{R}%
\right\Vert _{1}$\textquotedblright\ of $\mathbf{R}$ is defined as $\left\Vert
\mathbf{R}\right\Vert _{1}=\sum_{i,j=1}^{m}\left\vert \rho_{ij}\right\vert $.
We remark that, even though the assertion (\ref{eqSLLNV})\ is not explicitly
stated by Theorem 1 of \cite{Chen:2014SLLN}, it is written in the proof of
this theorem.

As a SLLN is perhaps the strongest characterization of the stability of a
sequence random variables, in this work we continue the line of research of
\cite{Chen:2014SLLN}, and characterize the type of dependence (via the order
of $\left\Vert \mathbf{R}\right\Vert _{1}$) under which (\ref{eqSLLNR}) and
(\ref{eqSLLNV}) hold when $\left(  \zeta_{i},\zeta_{j}\right)  ,i\neq j$
follows a Lancaster (but non-Gaussian) bivariate distribution with an infinite
support. It turns out that the strategy of \cite{Chen:2014SLLN} applies to the
settings here. Specifically, to prove (\ref{eqSLLNR}) we only need to
implement the following two steps: first, obtain a \textquotedblleft
comparison inequality\textquotedblright, i.e.,%
\begin{equation}
\left\vert \mathrm{cov}\left(  1\left\{  p_{i}\leq t\right\}  ,1\left\{
p_{j}\leq t\right\}  \right)  \right\vert \leq C\left\vert \rho_{ij}%
\right\vert \text{ for all }i\neq j\text{ and a constant }C>0; \label{IneqA}%
\end{equation}
second, apply Theorem 1 of \cite{Lyons:1988}, under the condition (\ref{eq1}),
to the indicators $X_{i}=1\left\{  p_{i}\leq t\right\}  $ with $1\leq i\leq m$
(or $i\in I_{0,m}$) that induces $R_{m}\left(  t\right)  $ (or $V_{m}\left(
t\right)  $). Once (\ref{eqSLLNR}) is proved and $\liminf_{m\rightarrow\infty
}m_{0}m^{-1}>0$ holds, (\ref{eqSLLNV}) follows as an easy corollary.

Our main result is the following:

\begin{theorem}
\label{ThmMain}Suppose that each pair $\left(  \zeta_{i},\zeta_{j}\right)
,i\neq j$ follows any of the following four Lancaster bivariate distributions
with correlation $\rho_{ij}$:

\begin{enumerate}
\item A Lancaster bivariate gamma distribution (defined by (\ref{defbiG}))
with shape parameter $\alpha\in(0,1]$;

\item A Lancaster bivariate Poisson distribution (defined by (\ref{biPoisson}%
)) with parameter $\alpha>0$;

\item A Lancaster negative binomial distribution (defined by (\ref{biNB}))
with parameter $\left(  \beta,c\right)  $ such that $\beta>0$ and $0<c<1$;

\item A Lancaster bivariate gamma-negative binomial distribution (defined by
(\ref{biGNB})) with parameter $\left(  \beta,c,\alpha\right)  $ such that
$\beta>0$, $0<c<1$ and $\alpha>0$.
\end{enumerate}

Then (\ref{IneqA}) holds. If (\ref{eq1}) holds, then (\ref{eqSLLNR}) holds. If
in addition $\liminf_{m\rightarrow\infty}m_{0}m^{-1}>0$, then (\ref{eqSLLNV}) holds.
\end{theorem}

The definitions of the four Lancaster bivariate distributions covered by
\autoref{ThmMain} can be found in \cite{Koudou:1998} and will be provided in
the proof of this theorem. Our findings seem to suggest that the inequality
(\ref{IneqA}) is universal for Lancaster bivariate distributions with infinite
supports. On the other hand, the Lancaster distributions considered by
\autoref{ThmMain} are often associated with the true null hypotheses in a
multiple testing scenario. For example, the Lancaster bivariate gamma
distribution includes the Lancaster bivariate central chi-square distribution
as a special case, and the latter distribution corresponds to the true null
hypothesis that its two marginal distributions have a zero centrality
parameter and the same degrees of freedom. Further, bivariate Poisson or
bivariate negative binomial distributions are widely used to model count data,
and the Lancaster bivariate Poisson or negative binomial distribution
corresponds to the true null hypothesis that its two marginal distributions
have identical parameters.

In view of the above discussion, \autoref{ThmMain} has the following
implication. Consider the slightly extended multiple testing scenario, where

\begin{itemize}
\item There are $\tilde{m}$ ($\geq m$) null hypotheses, $H_{i0}:\mu_{i}%
=\mu_{0}$ with $1\leq i\leq\tilde{m}$, to test simultaneously, each of which
has an associated test statistic $\zeta_{i}$;

\item Each $H_{i0}$ with $1\leq i\leq m$ is a true null hypothesis, and the
rest $\tilde{m}-m$ null hypotheses are false;

\item The MTP rejects an $H_{i0}$ iff its associated p-value $p_{i}\leq t$ for
a fixed rejection threshold $t\in\left(  0,1\right)  $.
\end{itemize}

\noindent Note that the above arrangement of the indices for the true and
false null hypotheses is unrestrictive. In this setting, the number of false
rejections of the MTP is $V_{\tilde{m}}\left(  t\right)  =\sum_{i=1}%
^{\tilde{m}}1\left\{  p_{i}\leq t\right\}  $, and the FDP of the MTP is%
\[
\mathrm{FDP}_{\tilde{m}}\left(  t\right)  =\frac{V_{\tilde{m}}\left(
t\right)  }{R_{\tilde{m}}\left(  t\right)  \vee1}\text{ \ with \ }R_{\tilde
{m}}\left(  t\right)  =\sum_{i=1}^{\tilde{m}}1\left\{  p_{i}\leq t\right\}  .
\]
Let $\mathbf{S}$ be the correlation matrix of $\left\{  \zeta_{i}\right\}
_{i=1}^{\tilde{m}}$ and $\pi_{0,\tilde{m}}=m^{-1}\tilde{m}$. When
$\liminf_{m\rightarrow\infty}\pi_{0,\tilde{m}}>0$,%
\begin{equation}
\tilde{m}^{-2}\left\Vert \mathbf{S}\right\Vert _{1}=O\big(\tilde{m}^{-\delta
}\big)\text{ \ for \ some}\ \delta>0\label{eq2a}%
\end{equation}
and the $p_{i}$'s associated with $H_{i0}$ for $1\leq i\leq m$ are identically
distributed as $p_{0}$, \autoref{ThmMain} implies
\begin{equation}
\lim_{m\rightarrow\infty}\left\vert m^{-1}V_{\tilde{m}}\left(  t\right)
-\mathbb{P}\left(  \left\{  p_{0}\leq t\right\}  \right)  \right\vert
=0\text{ almost surely.}\label{eq3}%
\end{equation}
Let $\hat{\pi}_{0,\tilde{m}}$ be an estimator of $\pi_{0,\tilde{m}}$ and set
\begin{equation}
\vartheta_{\tilde{m}}\left(  t\right)  =\frac{\hat{\pi}_{0,\tilde{m}%
}\mathbb{P}\left(  \left\{  p_{0}\leq t\right\}  \right)  }{\tilde{m}%
^{-1}R_{\tilde{m}}\left(  t\right)  }.\label{eq4}%
\end{equation}
It is easy to verify that, if $\hat{\pi}_{0,\tilde{m}}\pi_{0,\tilde{m}}%
^{-1}\rightsquigarrow1$ as $\tilde{m}\rightarrow\infty$, $\liminf_{\tilde
{m}\rightarrow\infty}\tilde{m}^{-1}R_{\tilde{m}}\left(  t\right)  >0$ almost surely and
(\ref{eq2a}) holds, then\ $\left\vert \vartheta_{\tilde{m}}\left(  t\right)
-\mathrm{FDP}_{\tilde{m}}\left(  t\right)  \right\vert \rightsquigarrow0$ as
$\tilde{m}\rightarrow\infty$, where $\rightsquigarrow$ denotes
\textquotedblleft convergence in probability\textquotedblright. Namely,
$\vartheta_{\tilde{m}}\left(  t\right)  $ consistently estimates
$\mathrm{FDP}_{\tilde{m}}\left(  t\right)  $ for each fixed $t\in\left(
0,1\right)  $. Note that $\vartheta_{\tilde{m}}\left(  t\right)  $ in
(\ref{eq4}) can be regarded as a slight extension of the FDR\ estimator
proposed by \cite{Storey:2004}.

A second implication of \autoref{ThmMain} is as follows. The \textquotedblleft
weak dependence" assumption, proposed in \cite{Storey:2004} and widely used in
the multiple testing literature, requires that there exist two continuous
functions $G_{0}$ and $G_{1}$ such that\ for each $t\in(0,1],$%
\begin{equation}
\lim_{m\rightarrow\infty}m_{0}^{-1}V_{m}\left(  t\right)  =G_{0}\left(
t\right)  \text{ \ and \ }\lim_{m-\infty}\left(  m-m_{0}\right)  ^{-1}\left[
R_{m}\left(  t\right)  -V_{m}\left(  t\right)  \right]  =G_{1}\left(
t\right)  \label{eq5}%
\end{equation}
almost surely. However, to check whether (\ref{eq5}) holds is often very hard
(even after the continuity requirement on $G_{0}$ and $G_{1}$ is removed).
\autoref{ThmMain} here and Theorem 1 in \cite{Chen:2014SLLN} together provide
a way to check whether this assumption holds in the scenario of simultaneously
testing the parameters of a larger number of dependent random variables, each
pair of which follows any of the five Lancaster bivariate distributions that
are studied in \cite{Koudou:1998}. Specifically, a check on the order of the
$l_{1}$-norm of the correlation matrix of these random variables suffices for
this purpose. We will report in another article on how to consistently
estimate $m^{-2}\left\Vert \mathbf{R}\right\Vert _{1}$ or efficiently test the
order of $\left\Vert \mathbf{R}\right\Vert _{1}$.

The rest of the article is devoted to the proof of \autoref{ThmMain}.

\section{Proof of \autoref{ThmMain}}

In the proof, $\mathbb{V}\left[  \cdot\right]  $ and $\mathsf{cov}\left[
\cdot,\cdot\right]  $ are the variance and covariance operators,
$\mathbb{N}_{0}=\mathbb{N}\cup\left\{  0\right\}  $, and $C$ denotes a
positive constant that can assume different (and appropriate) values at
different occurrences. We need Theorem 1 of \cite{Lyons:1988} in the proof,
which reads \textquotedblleft Let $\left\{  \chi_{n}\right\}  _{n=1}^{\infty}$
be a sequence of complex-valued random variables such that $\mathbb{E}\left[
\left\vert \chi_{n}\right\vert ^{2}\right]  \leq1$. Set \thinspace
$Q_{N}=N^{-1}\sum\nolimits_{n=1}^{N}\chi_{n}$. If $\left\vert \chi
_{n}\right\vert \leq1$ a.s. and%
\begin{equation}
\sum\nolimits_{N=1}^{\infty}N^{-1}\mathbb{E}\left[  \left\vert Q_{N}%
\right\vert ^{2}\right]  <\infty, \label{eqCondLyons}%
\end{equation}
then $\lim_{N\rightarrow\infty}Q_{N}=0$ a.s.\textquotedblright\ \noindent A
sufficient condition for the SLLN to hold for $\left\{  \chi_{n}\right\}
_{n=1}^{\infty}$ is that $\mathbb{E}\left[  \left\vert Q_{m}\right\vert
^{2}\right]  =O\left(  m^{-\delta}\right)  $ for some $\delta>0$, which
implies (\ref{eqCondLyons}).

Now we present the arguments. Recall $X_{i}=1\left\{  p_{i}\leq t\right\}  $,
for which $R_{m}\left(  t\right)  =\sum_{i=1}^{m}X_{i}$ and $V_{m}\left(
t\right)  =\sum_{i\in I_{0,m}}X_{i}$. We aim to show that $\mathbb{V}\left[
m^{-1}R_{m}\left(  t\right)  \right]  $ satisfies $O\left(  m^{-\delta_{\ast}%
}\right)  $ with $\delta_{\ast}=\min\left\{  \delta,1\right\}  $. Define two
sets%
\[
E_{1,m}=\left\{  \left(  i,j\right)  :1\leq i,j\leq m,i\neq j,\left\vert
\rho_{ij}\right\vert =1\right\}
\]
and%
\[
E_{2,m}=\left\{  \left(  i,j\right)  :1\leq i,j\leq m,i\neq j,\left\vert
\rho_{ij}\right\vert <1\right\}  .
\]
Namely, $E_{2,m}$ records pairs $\left(  \zeta_{i},\zeta_{j}\right)  $ with
$i\neq j$ such that $\zeta_{i}$ and $\zeta_{j}$ are linearly dependent almost
surely. Obviously, $\left\vert \mathrm{cov}\left(  X_{i},X_{j}\right)
\right\vert \leq C=C\left\vert \rho_{ij}\right\vert $ for $\left(  i,j\right)
\in E_{2,m}$. Further,
\begin{align}
\mathbb{V}\left[  m^{-1}R_{m}\left(  t\right)  \right]   &  \leq O\left(
m^{-1}\right)  +m^{-2}\sum_{\left(  i,j\right)  \in E_{2,m}}\left\vert
\mathrm{cov}\left(  X_{i},X_{j}\right)  \right\vert +m^{-2}\sum_{\left(
i,j\right)  \in E_{1,m}}\left\vert \mathrm{cov}\left(  X_{i},X_{j}\right)
\right\vert \nonumber\\
&  \leq O\left(  m^{-\min\left\{  \delta,1\right\}  }\right)  +m^{-2}%
\sum_{\left(  i,j\right)  \in E_{1,m}}\left\vert \mathrm{cov}\left(
X_{i},X_{j}\right)  \right\vert \label{eqE1}%
\end{align}
since%
\[
m^{-2}\sum\nolimits_{\left(  i,j\right)  \in E_{2,m}}\left\vert \mathrm{cov}%
\left(  X_{i},X_{j}\right)  \right\vert =O\left(  m^{-2}\left\Vert
\mathbf{R}\right\Vert _{1}\right)  =O\big(m^{-\delta}\big).
\]
So, we only need to upper bound $B_{1,m}=m^{-2}\sum_{\left(  i,j\right)  \in
E_{1,m}}\left\vert \mathrm{cov}\left(  X_{i},X_{j}\right)  \right\vert $ on
the right-hand side of (\ref{eqE1}).

On the other hand,
\begin{align*}
\mathbb{V}\left[  m^{-1}V_{m}\left(  t\right)  \right]   &  \leq O\left(
m^{-1}\right)  +m^{-2}\sum_{\left(  i,j\right)  \in\tilde{E}_{2,m}}\left\vert
\mathrm{cov}\left(  X_{i},X_{j}\right)  \right\vert +m^{-2}\sum_{\left(
i,j\right)  \in\tilde{E}_{1,m}}\left\vert \mathrm{cov}\left(  X_{i}%
,X_{j}\right)  \right\vert \\
&  \leq Cm^{-1}+Cm^{-\delta}+CB_{1,m},
\end{align*}
where $\tilde{E}_{k,m}=E_{k,m}\cap\left(  I_{0,m}\times I_{0,m}\right)  $ for
$k\in\left\{  1,2\right\}  $. So, an upper bound on $B_{1,m}$ will induce the
same upper bound for $\mathbb{V}\left[  m^{-1}R_{m}\left(  t\right)  \right]
$ and $\mathbb{V}\left[  m^{-1}V_{m}\left(  t\right)  \right]  $.

We will split the rest of the proof into four cases in terms of upper bounding
$B_{1,m}$, each corresponding to a Lancaster bivariate distribution in the
statement of \autoref{ThmMain} and each occupying a subsection.

\subsection{The Lancaster bivariate gamma distribution}

\label{lancaster_bvgamma}

The Lancaster bivariate gamma distribution was derived by \cite{griffiths1969}
and \cite{Koudou:1998}. Specifically, if $\left(  X,Y\right)  $ follows this
distribution with shape parameter $\alpha>0$ and correlation $\rho\in
\lbrack0,1)$, then its density is%
\begin{equation}
h\left(  x,y;\alpha,\rho\right)  =f\left(  x;\alpha\right)  f\left(
y;\alpha\right)  \sum_{n=0}^{\infty}\frac{\rho^{n}n!}{\Gamma\left(
\alpha+n\right)  \Gamma\left(  \alpha\right)  }L_{n}^{\left(  \alpha-1\right)
}\left(  x\right)  L_{n}^{\left(  \alpha-1\right)  }\left(  y\right)  ,
\label{defbiG}%
\end{equation}
where%
\[
f\left(  x;\alpha\right)  =\frac{1}{\Gamma\left(  \alpha\right)  }x^{\alpha
-1}e^{-x}\text{ \ for \ }x>0
\]
is the gamma density with shape parameter $\alpha>0$, and
\[
L_{n}^{\left(  \alpha\right)  }\left(  x\right)  =\sum_{k=0}^{n}%
\binom{n+\alpha}{n-k}\frac{\left(  -x\right)  ^{k}}{k!}\text{ \ for \ }%
n\in\mathbb{N}_{0}%
\]
is the $n$th Laguerre polynomial of order $\alpha>0$.

Let $\tau=F_{i}^{-1}\left(  1-t\right)  $. If $\left(  \zeta_{i},\zeta
_{j}\right)  $ with $\left(  i,j\right)  \in E_{2,m}$ follows a Lancaster
bivariate gamma distribution with shape parameter $\alpha>0$ and correlation
$\rho_{ij}\in\lbrack0,1)$, then%
\[
\kappa_{ij}=\mathrm{cov}\left(  1\left\{  p_{i}\leq t\right\}  ,1\left\{
p_{j}\leq t\right\}  \right)  =\sum_{n=1}^{\infty}\frac{\rho_{ij}^{n}%
n!}{\Gamma\left(  \alpha+n\right)  \Gamma\left(  \alpha\right)  }q_{n}%
^{2}\left(  \tau;\alpha\right)  ,
\]
where%
\[
q_{n}\left(  \tau;\alpha\right)  =\int_{-\infty}^{\tau}f\left(  x;\alpha
\right)  L_{n}^{\left(  \alpha-1\right)  }\left(  x\right)  dx=\frac{1}%
{\Gamma\left(  \alpha\right)  }\int_{-\infty}^{\tau}x^{\alpha-1}e^{-x}%
L_{n}^{\left(  \alpha-1\right)  }\left(  x\right)  dx.
\]
From the Rodrigue's formula (e.g., on page 101 of \cite{Szego:1939}), i.e.,%
\[
L_{n}^{(\alpha)}(x)=\frac{1}{n!}x^{-\alpha}e^{x}\frac{d^{n}}{dx^{n}}\left(
x^{n+\alpha}e^{-x}\right)  \text{ \ \ for }n\in\mathbb{N}_{0}\text{ \ and
}\alpha>-1,
\]
we obtain, for $y>0$ and $n\geq1$,%
\begin{align*}
\int_{-\infty}^{y}x^{\alpha}e^{-x}L_{n}^{(\alpha)}(x)dx &  =\frac{1}{n!}%
\int_{-\infty}^{y}\left[  \frac{d^{n}}{dx^{n}}\left(  x^{n+\alpha}%
e^{-x}\right)  \right]  dx\\
&  =\frac{y^{\alpha+1}e^{-y}}{n}\frac{y^{-\left(  \alpha+1\right)  }e^{y}%
}{\left(  n-1\right)  !}\left[  \left.  \frac{d^{n-1}}{dx^{n-1}}\left(
x^{n-1+\alpha+1}e^{-x}\right)  \right\vert _{x=y}\right]  \\
&  =\frac{y^{\alpha+1}e^{-y}}{n}L_{n-1}^{(\alpha+1)}(y).
\end{align*}
Therefore,%
\[
\kappa_{ij}=\sum_{n=1}^{\infty}\frac{\rho_{ij}^{n}n!}{n^{2}\Gamma\left(
\alpha+n\right)  \Gamma^{3}\left(  \alpha\right)  }\left[  \tau^{\alpha
}e^{-\tau}L_{n-1}^{(\alpha)}(\tau)\right]  ^{2}.
\]
By Watson's bound on page 21 of \cite{Watson1939jlms}, i.e.,%
\[
\left\vert L_{n}^{\left(  \alpha\right)  }\left(  x\right)  \right\vert
\leq\frac{\Gamma\left(  \alpha+1+n\right)  }{\Gamma\left(  \alpha+1\right)
n!}e^{x/2}\text{ \ for \ }x\geq0,\alpha\geq0\text{ and }n\in\mathbb{N}_{0},
\]
we obtain%
\[
\left\vert \kappa_{ij}\right\vert \leq C\sum_{n=1}^{\infty}\frac{\rho_{ij}%
^{n}\Gamma\left(  \alpha+n\right)  }{n!}\tau^{2\alpha}e^{-\tau}.
\]
However, the identity (1) in \cite{tricomi1951asymptotic} states that,\ for
distinct real constants $\alpha$ and $\gamma,$
\begin{equation}
\frac{\Gamma\left(  z+\alpha\right)  }{\Gamma\left(  z+\gamma\right)
}=z^{\alpha-\gamma}\left[  1+\frac{\left(  \alpha-\gamma\right)  \left(
\alpha+\gamma-1\right)  }{2z}+O\left(  \left\vert z\right\vert ^{-2}\right)
\right]  \text{ \ as \ }\left\vert z\right\vert \rightarrow\infty.\label{eq8}%
\end{equation}
So, when $\alpha\leq1$,%
\[
\left\vert \kappa_{ij}\right\vert \leq C\sum_{n=1}^{\infty}\frac{\rho_{ij}%
^{n}}{n^{1-\alpha}}=C\rho_{ij}\sum_{n=1}^{\infty}\frac{\rho_{ij}^{n-1}%
}{n^{1-\alpha}}\leq C\rho_{ij},
\]
and (\ref{eqSLLNR}) holds.

\begin{remark}
If $\zeta_{i}$ is the central chi-square random variable with $v$ degrees of
freedom and density%
\[
f\left(  x;v/2\right)  =\frac{1}{\Gamma\left(  v/2\right)  2^{v/2}}%
x^{v/2-1}e^{-x/2}\text{ \ for \ }x>0,
\]
then \autoref{ThmMain} is valid when $v=1$ or $2$.
\end{remark}

\subsection{The Lancaster bivariate Poisson distribution}

\label{lancaster_bvpoisson}

For $a>0$ and $x,n\in\mathbb{N}_{0}$, let
\[
C_{n}(x;a)=\sqrt{\frac{a^{n}}{n!}}\sum_{k=0}^{n}(-1)^{k}\binom{n}{k}\binom
{x}{k}\,\frac{k!}{a^{k}}%
\]
denote the Charlier polynomial of degree $n$, where $\binom{n}{k}=\frac
{n!}{k!\left(  n-k\right)  !}$ if $n\geq k$ and $\binom{n}{k}=0$ if $n<k$. The
Lancaster bivariate Poisson distribution was derived by \cite{Koudou:1998}.
Specifically, if $\left(  X,Y\right)  $ follows such a distribution with
parameter $a>0$ and corelation $\lambda\in\lbrack0,1]$, then it has density
\begin{equation}
h\left(  x,y;a,\rho\right)  =f(x;a)\,f(y;a)\sum_{n=0}^{\infty}\rho^{n}%
\,C_{n}(x;a)\,C_{n}(y;a)\text{ \ for }x,y\in\mathbb{N}_{0}, \label{biPoisson}%
\end{equation}
where
\[
f(x;a)=\frac{a^{x}\,e^{-a}}{x!}\text{ for }x\in\mathbb{N}_{0}%
\]
is the \textit{probability mass function (PMF)} for a Poisson random variable
with mean $a$.

Set $\tau=F_{i}^{-1}\left(  1-t\right)  $, and let $x_{0}$ be the integer part
of $\tau$. If $\left(  \zeta_{i},\zeta_{j}\right)  \ $with $\left(
i,j\right)  \in E_{2,m}$ follows a Lancaster bivariate Poisson distribution
with correlation $\rho_{ij}\in\left[  0,1\right]  $, then%
\[
\kappa_{ij}=\mathrm{cov}\left(  1\left\{  p_{i}\leq t\right\}  ,1\left\{
p_{j}\leq t\right\}  \right)  =\sum_{n=1}^{\infty}\rho_{ij}^{n}q_{n}%
^{2}\left(  x_{0};a\right)  ,
\]
where%
\[
q_{n}\left(  x_{0};a\right)  =\sum_{x=0}^{x_{0}}f(x;a)C_{n}(x;a)=\sqrt
{\frac{a^{n}}{n!}}\sum_{x=0}^{x_{0}}\frac{a^{x}\,e^{-a}}{x!}\sum_{k=0}%
^{n}(-1)^{k}\binom{n}{k}\binom{x}{k}\,\frac{k!}{a^{k}}.
\]
It suffices to bound $q_{n}\left(  x_{0};a\right)  $. Specifically,
\[
\left\vert q_{n}\left(  x_{0};a\right)  \right\vert \leq\sqrt{\frac{a^{n}}%
{n!}}\sum_{x=0}^{x_{0}}a^{x}\,e^{-a}\sum_{k=0}^{n}\frac{a^{-k}}{\left(
x-k\right)  !}\binom{n}{k}\,\leq C\sqrt{\frac{a^{n}}{n!}}\left(
1+a^{-1}\right)  ^{n}.
\]
So,%
\[
\left\vert \kappa_{ij}\right\vert \leq C\sum_{n=1}^{\infty}\frac{a^{n}%
\rho_{ij}^{n}}{n!}\left(  1+a^{-1}\right)  ^{2n}\leq C\rho_{ij}%
\]
and (\ref{eqSLLNR}) holds.

\subsection{The Lancaster bivariate negative binomial distribution}

\label{lancaster_bvNegBin}

Let $\beta>0$ and $0<c<1$, and $M_{n}^{\beta,c}(x)$ denote the $nt$h
(normalized) Meixner polynomial, i.e.,
\[
M_{n}^{\beta,c}(x)=\sqrt{\frac{c^{n}\,(\beta)_{n}}{n!}}\sum_{k=0}^{n}%
\frac{(-n)_{k}\,(-x)_{k}}{(\beta)_{k}\,k!}\left(  1-c^{-1}\right)  ^{k}\text{
\ for \ }x\in\mathbb{N}_{0}.
\]
Here $\left(  a\right)  _{n}=\prod_{k=0}^{n-1}\left(  a+k\right)  $ for
$a\in\mathbb{R}$ and $n\in\mathbb{N}$, and $\left(  -x\right)  _{k}=0$ is set
when $x<k$. The Lancaster bivariate negative binomial distribution with
parameter $\left(  \beta,c\right)  $ and correlation $\rho\in\lbrack0,1)$ was
derived by \cite{Koudou:1998}. Specifically, if $(X,Y)$ follows such a
distribution, then it has density
\begin{equation}
h\left(  x,y;\beta,c\right)  =f(x;\beta,c)\,f(y;\beta,c)\sum_{n=0}^{\infty
}\rho^{n}\,M_{n}^{\beta,c}(x)\,M_{n}^{\beta,c}(y)\ \ \text{for }%
x,y\in\mathbb{N}_{0}\text{ \ and \ }0\leq\rho<1, \label{biNB}%
\end{equation}
where
\[
f(x;\beta,c)=(1-c)^{\beta}\,\frac{c^{x}\,(\beta)_{x}}{x!}\text{ \ for }%
x\in\mathbb{N}_{0}%
\]
is the PMF for a negative binomial random variable.

Set $\tau=F_{i}^{-1}\left(  1-t\right)  $, and let $x_{0}$ be the integer part
of $\tau$. If $\left(  \zeta_{i},\zeta_{j}\right)  \ $with $\left(
i,j\right)  \in E_{2,m}$ follows a Lancaster bivariate negative binomial
distribution with parameter $\left(  \beta,c\right)  $ and correlation
$\rho_{ij}\in\lbrack0,1)$, then%
\[
\kappa_{ij}=\mathrm{cov}\left(  1\left\{  p_{i}\leq t\right\}  ,1\left\{
p_{j}\leq t\right\}  \right)  =\sum_{n=1}^{\infty}\rho_{ij}^{n}q_{n}%
^{2}\left(  x_{0};\beta,c\right)  ,
\]
where%
\begin{align*}
q_{n}\left(  x_{0};\beta,c\right)   &  =\sum_{x=0}^{x_{0}}f(x;\beta
,c)M_{n}^{\beta,c}(x)\\
&  =\sqrt{\frac{c^{n}\,(\beta)_{n}}{n!}}\sum_{x=0}^{x_{0}}(1-c)^{\beta}%
\,\frac{c^{x}\,(\beta)_{x}}{x!}\sum_{k=0}^{n}\frac{(-n)_{k}\,(-x)_{k}}%
{(\beta)_{k}\,k!}\left(  1-c^{-1}\right)  ^{k}.
\end{align*}
It suffices to bound $q_{n}\left(  x_{0};\beta,c\right)  $. Specifically,
\begin{align*}
\left\vert q_{n}\left(  x_{0};\beta,c\right)  \right\vert  &  \leq
C\sqrt{\frac{c^{n}\,(\beta)_{n}}{n!}}\sum_{k=0}^{x_{0}}\frac{n\left(
n-1\right)  \cdots\left(  n-k+1\right)  \,}{(\beta)_{k}\,k!}\left\vert
1-c^{-1}\right\vert ^{k}\\
&  \leq C\sqrt{\frac{c^{n}\,(\beta)_{n}}{n!}}x_{0}n^{x_{0}}\leq Cc^{n/2}%
n^{\left(  \beta-1+2x_{0}\right)  /2},
\end{align*}
where we have applied the identity (\ref{eq8}) to obtain the last inequality.
Since $0<c<1$, we have%
\[
\left\vert \kappa_{ij}\right\vert \leq C\sum_{n=1}^{\infty}\rho_{ij}^{n}%
c^{n}n^{\beta-1+2x_{0}}\leq C\rho_{ij}.
\]
So, (\ref{eqSLLNR}) holds.

\subsection{The Lancaster bivariate gamma-negative binomial distribution}

Let $\sigma>0$, $\beta>0$ and $0<c<1$ be three constants. The Lancaster
bivariate gamma-negative binomial distribution was derived by
\cite{Koudou:1998}. Specifically, if $(X,Y)$ follows such a distribution with
parameter $\left(  \alpha,\beta,c\right)  $ and correlation $\rho\in\left[
0,\sqrt{c}\right]  $, then it has density
\begin{equation}
h\left(  x,y;\alpha,\beta,c\right)  =f(x;\beta,c)\,g\left(  y;\alpha\right)
\sum_{n=0}^{\infty}\rho^{n}\,\sqrt{\frac{n!}{\left(  \alpha\right)  _{n}}%
}L_{n}^{\left(  \alpha-1\right)  }\left(  x\right)  \,M_{n}^{\beta,c}(y)
\label{biGNB}%
\end{equation}
for $x\in\mathbb{N}_{0}$ and $y>0$, for which $X$ is a negative binomial
random variable with PMF
\[
f(x;\beta,c)=(1-c)^{\beta}\,\frac{c^{x}\,(\beta)_{x}}{x!}\text{ \ for \ }%
x\in\mathbb{N}_{0},
\]
and $Y$ is a gamma random variable with density
\[
g\left(  y;\alpha\right)  =\frac{1}{\Gamma\left(  \alpha\right)  }y^{\alpha
-1}e^{-y}\text{ \ for \ }y>0.
\]
\

If $\left(  \zeta_{i},\zeta_{j}\right)  \ $with $\left(  i,j\right)  \in
E_{2,m}$ follows a Lancaster bivariate gamma-negative binomial distribution
with parameter $\left(  \alpha,\beta,c\right)  $ and correlation $\rho_{ij}%
\in\lbrack0,1)$, then%
\[
\kappa_{ij}=\mathrm{cov}\left(  1\left\{  p_{i}\leq t\right\}  ,1\left\{
p_{j}\leq t\right\}  \right)  =\sum_{n=1}^{\infty}\rho_{ij}^{n}q_{n}\left(
x_{0};\beta,c\right)  r_{n}\left(  \tau,\alpha\right)  ,
\]
where $x_{0}$ is the integer part of $F_{i}^{-1}\left(  1-t\right)  $,
$\tau=F_{j}^{-1}\left(  1-t\right)  $,%
\[
q_{n}\left(  x_{0};\beta,c\right)  =\sum_{x=0}^{x_{0}}f(x;\beta,c)M_{n}%
^{\beta,c}(x)
\]
and%
\[
r_{n}\left(  \tau,\alpha\right)  =\sqrt{\frac{n!}{\left(  \alpha\right)  _{n}%
}}\int_{-\infty}^{\tau}g\left(  y;\alpha\right)  L_{n}^{\left(  \alpha
-1\right)  }\left(  x\right)  dx.
\]
Using the bounds derived in \autoref{lancaster_bvgamma} and
\autoref{lancaster_bvNegBin}, we obtain
\[
\left\vert \kappa_{ij}\right\vert \leq\sum_{n=1}^{\infty}\rho_{ij}^{n}%
c^{n/2}n^{\left(  \beta-1+2x_{0}\right)  /2}n^{\left(  1-\alpha\right)
/2}\leq C\rho_{ij}.
\]
So, (\ref{eqSLLNR}) holds.

\section*{Acknowledgements}

This research was funded by the New Faculty Seed Grant provided by Washington
State University. I would like to thank G\'{e}rard Letac for suggesting to me
several references on orthogonal polynomials, and Donald Richards for
encouraging me to extend the work of \cite{Chen:2014SLLN} to the setting of
Lancaster bivariate distributions.

\bibliographystyle{chicago}

\end{document}